\documentclass[11pt,a4paper]{article}
\usepackage[T1]{fontenc}
\usepackage{lmodern}
\usepackage[margin=28mm]{geometry}
\usepackage{amsmath,amssymb,amsthm}
\usepackage{mathrsfs}
\usepackage{enumitem}
\usepackage{needspace}
\usepackage{microtype}
\usepackage[colorlinks=true,linkcolor=blue,citecolor=blue,urlcolor=blue]{hyperref}
\hypersetup{pdftitle={Relative Rouquier linearity and inverse hybrid Kazhdan--Lusztig positivity},pdfauthor={Kaixuan Gan, Xun Xie, Qifan Zhou}}
\numberwithin{equation}{section}
\newtheorem{theorem}{Theorem}[section]
\newtheorem{proposition}[theorem]{Proposition}
\newtheorem{lemma}[theorem]{Lemma}

\theoremstyle{remark}

\newcommand{\Z}{\mathbb Z}
\newcommand{\N}{\mathbb N}
\newcommand{\B}{\mathcal B}
\newcommand{\HH}{\mathcal H}
\newcommand{\LL}{\mathcal L}
\newcommand{\ch}{\operatorname{ch}}
\newcommand{\Cone}{\operatorname{Cone}}
\newcommand{\summand}{\mathrel{\subseteq_{\oplus}}}
\setlist[enumerate,1]{label=\textup{(\roman*)},leftmargin=2em,itemsep=3pt,topsep=5pt}
\title{Relative Rouquier linearity and inverse hybrid\\Kazhdan--Lusztig positivity}
\author{Kaixuan Gan, Xun Xie and Qifan Zhou%
\thanks{The authors are partially supported by NSFC Grant Nos.~12171030 and~12431002.}\\[0.5ex]
\normalsize Beijing Institute of Technology, Beijing, China}

\begin{document}
\maketitle
\vspace{-1.8em}
\begin{abstract}
We prove signed positivity for the expansion of a parabolic hybrid
basis in the Kazhdan--Lusztig basis of a Hecke algebra.
The coefficients have degree bounded by the length of
the corresponding minimal coset representative. 
We deduce this linearity from Elias--Williamson's ordinary Rouquier
linearity by a truncation argument. The result does not require the parabolic subgroup
to be finite.
\end{abstract}

\section{Introduction and main result}

The parabolic hybrid basis combines the Kazhdan--Lusztig basis of a
standard parabolic subgroup with standard basis elements indexed by
minimal coset representatives. We study its expansion in the
Kazhdan--Lusztig basis and show that the coefficients have signs
prescribed by differences of lengths.

Let $(W,S)$ be a Coxeter system with $S$ finite, identity $e$, length function
$\ell$, and Bruhat order $\le$. Let $\HH$ be the Hecke algebra over
$\Z[v,v^{-1}]$, with standard basis $\{H_w\}_{w\in W}$, normalized by
\[
 H_s^2=1+(v^{-1}-v)H_s.
\]
The bar involution is given by $\overline v=v^{-1}$ and
$\overline{H_w}=H_{w^{-1}}^{-1}$. The Kazhdan--Lusztig basis
$\{C_w\}_{w\in W}$ is characterized by
\[
 \overline{C_w}=C_w,\qquad
 C_w=H_w+\sum_{u<w}h_{u,w}(v)H_u,
 \qquad h_{u,w}(v)\in v\Z[v].
\]
In particular, $C_s=H_s+v$. We use $\N=\{0,1,2,\ldots\}$ and write
$[v^i]p$ for the coefficient of $v^i$ in a Laurent polynomial $p$.

For $I\subseteq S$, put
\[
 W_I=\langle I\rangle,\qquad
 {}^I W=\{y\in W:\ell(sy)>\ell(y)\text{ for all }s\in I\}.
\]
Thus ${}^I W$ is the set of minimal representatives for
$W_I\backslash W$. Every $w\in W$ has a unique length-additive factorization $w=xy$, where
$x\in W_I$ and $y\in{}^I W$. The elements $C_xH_y$ form a basis,
called the parabolic hybrid basis following \cite[\S3.1]{GH}.
Positivity of the expansion of the KL basis in the hybrid
basis is proved in \cite[Corollary~1.2]{XZ}.
We study the inverse transition matrix, which expresses the hybrid
basis in the KL basis. Write
$r^I_{z;x,y}$ for its entries and $Q^I_{z;x,y}$ for their sign
normalizations:
\begin{equation}\label{eq:inverse}
 C_xH_y=\sum_{z\in W}r^I_{z;x,y}(v)C_z,
 \qquad
 Q^I_{z;x,y}(v)=(-1)^{\ell(xy)-\ell(z)}r^I_{z;x,y}(v).
\end{equation}

\Needspace{10\baselineskip}
\begin{theorem}[Inverse hybrid positivity]\label{thm:positive}
Let $I\subseteq S$, $x\in W_I$, and $y\in{}^I W$.
\begin{enumerate}
 \item $Q^I_{xy;x,y}=1$, and $Q^I_{z;x,y}=0$ if $z\not\le xy$.
 \item If $z<xy$, then $Q^I_{z;x,y}\in v\N[v]$.
 \item If $[v^i]Q^I_{z;x,y}\ne0$, then
 \[
  0\le i\le\ell(y),\qquad
  i\equiv\ell(xy)-\ell(z)\pmod 2.
 \]
\end{enumerate}
\end{theorem}

We give a short proof of this result using ordinary Rouquier
linearity. For $I=\varnothing$, it reduces to ordinary inverse
Kazhdan--Lusztig positivity \cite[Remark 6.10]{EW}.
If $W_I$ is finite and $x=w_I$ is its longest element, the
anti-involution $H_w\mapsto H_{w^{-1}}$ sends $C_{w_I}H_y$ to
$H_{y^{-1}}C_{w_I}$. Thus the case $x=w_I$ recovers the
positivity of inverse parabolic Kazhdan--Lusztig polynomials
\cite[Theorem 4.15 and Remark 4.16]{P}.

The proof reduces positivity to the linearity of the minimal
model of $B_xF_y$, where $B_x$ is an indecomposable Soergel
bimodule and $F_y$ is a Rouquier complex.
In Section~\ref{sec:linearity}, we prove this linearity by
comparing internal shifts with cohomological degrees.
An inequality from \cite[Lemma 6.6]{EW} gives one bound,
while a truncation of $F_x$, together with induction on
$\ell(x)$, gives the other.
In Section~\ref{sec:positivity}, we use Euler characters
and parity to identify the coefficients of the sign-normalized
polynomials with the multiplicities of indecomposable summands.

\section{Relative linearity}\label{sec:linearity}

Choose the real reflection faithful realization $V$ of \cite[\S3.1]{EW}, and let $\B$ be the
corresponding category of Soergel bimodules. Set
$R=\operatorname{Sym}(V^*)$, with linear functions in degree $2$.
The internal shift of a graded vector space $M$ is written $M(a)$, with $M(a)_d=M_{d+a}$.
Let $B_w$ denote the self-dual indecomposable Soergel
bimodule indexed by $w\in W$. In particular, $B_e=R$ and
$B_s=R\otimes_{R^s}R(1)$, where $R^s$ is the subring of $s$-invariants.

We work in \(K^b(\mathcal B)\), the homotopy category of bounded complexes in \(\mathcal B\) with differentials of internal degree zero. For such a complex \(E\), we write \(E^i\) for its term in cohomological degree \(i\). We use the shift convention
\[
(E[k])^i=E^{i+k},\qquad d_{E[k]}=(-1)^k d_E.
\]
Tensor products are over $R$, and $\simeq$ denotes homotopy equivalence.
The Soergel character theorem \cite[Theorem 1.1]{EW} gives
$\ch(B_w(a))=v^aC_w$. We extend the character to bounded complexes by
\[
\operatorname{ch}(E)=\sum_i(-1)^i\operatorname{ch}(E^i).
\]

A bounded complex is called \textit{minimal} if it has no nonzero contractible direct summand. Since \(\mathcal B\) is Krull–Schmidt, this is equivalent to requiring that no differential has an isomorphism component between indecomposable direct summands. Moreover, every bounded complex is homotopy equivalent to a minimal complex, unique up to
isomorphism of complexes, which can be chosen as a direct summand of
the original complex \cite[\S6.1]{EW}. We write \(M\summand N\) to mean that \(M\) is isomorphic to a direct summand of \(N\), either in \(\mathcal B\) or in the category of bounded complexes, as appropriate.

A complex is \emph{linear} if its term in degree $i$ is a finite direct
sum of bimodules $B_z(i)$. The shift $(a)[-a]$ preserves linearity.

Let $F_s=[B_s\to R(1)]$, with $B_s$ in degree $0$ and multiplication
as differential. For a reduced expression $w=s_1\cdots s_d$, let
$F_w$ be a minimal model of $F_{s_1}\cdots F_{s_d}$, called the minimal Rouquier complex. Its isomorphism
class is independent of the reduced expression. If
$\ell(ab)=\ell(a)+\ell(b)$, then $F_aF_b\simeq F_{ab}$.

\Needspace{6\baselineskip}

By \cite[Theorem 6.9]{EW}, \(F_w\) satisfies
\[
F_w^0=B_w,\qquad
F_w^i\cong\bigoplus_{z<w}B_z(i)^{\oplus n_{z,i}}
\quad(i>0),
\]
where \(n_{z,i}\in\mathbb N\), and \(F_w^i=0\) for \(i<0\).
In particular, \(F_w\) is linear.

We also use \cite[Lemma 6.6]{EW}: if a minimal complex $E$ satisfies
\begin{equation}\label{eq:upper-condition}
 B_z(j)\summand E^i\quad\Longrightarrow\quad j\le i
\end{equation}
then the minimal model of $EF_s$ satisfies the same bound.

\begin{lemma}\label{lem:cones}
	The following statements hold in $K^b(\B)$.
	\begin{enumerate}
		\item If
		\[
		A\longrightarrow E\longrightarrow C
		\longrightarrow A[1]
		\]
		is a distinguished triangle and $A,C$ are homotopy
		equivalent to linear complexes, then so is $E$.
		
		\item Let $E$ be a minimal model of the cone of a
		chain map between linear complexes. Then
		\[
		B_z(j)\summand E^i
		\quad\Longrightarrow\quad
		j\in\{i,i+1\}.
		\]
		In particular, $j\ge i$.
	\end{enumerate}
\end{lemma}

\begin{proof}
	Recall that for a chain map $f:A\to B$,
	\(
	\operatorname{Cone}(f)^i=B^i\oplus A^{i+1}.
	\)
	
	(i) Replacing $A$ and $C$ by linear representatives,
	we may write
	\[
	E\simeq\operatorname{Cone}(C[-1]\to A).
	\]
	The term in degree $i$ of this cone is
	$A^i\oplus C^i$, which is a direct sum of
	bimodules $B_z(i)$. Thus $E$ is homotopy
	equivalent to a linear complex.
	
	(ii) Let $f:A\to B$ be a chain map between
	linear complexes. Every indecomposable summand
	of
	\[
	\operatorname{Cone}(f)^i=B^i\oplus A^{i+1}
	\]
	has internal shift $i$ or $i+1$.
	By the existence of minimal subcomplexes
	\cite[\S 6.1]{EW}, the terms of a minimal model
	are direct summands of the corresponding terms
	of $\operatorname{Cone}(f)$.
	The assertion follows.
\end{proof}

\begin{proposition}[Relative linearity]\label{prop:linear}
	Let $x\in W_I$ and $y\in{}^I W$, and let
	$\LL_{x,y}$ be a minimal model of $B_xF_y$.
	Then:
	\begin{enumerate}
		\item $\LL_{x,y}$ is linear. That is,
		\[
		\LL_{x,y}^i\cong
		\bigoplus_{z\in W}B_z(i)^{\oplus m_{z,i}},
		\qquad m_{z,i}\in\N,
		\]
		with only finitely many nonzero multiplicities.
		
		\item $\LL_{x,y}^i=0$ unless
		$0\le i\le\ell(y)$.
	\end{enumerate}
\end{proposition}

\begin{proof}
	(i) We prove that
	\[
	B_z(j)\summand\LL_{x,y}^i
	\quad\Longrightarrow\quad j=i.
	\]
	
	First, $B_x[0]$ satisfies the condition
	\eqref{eq:upper-condition}.
	Applying \cite[Lemma 6.6]{EW} successively
	along a reduced expression of $y$, we obtain
	\begin{equation}\label{eq:upper}
		B_z(j)\summand\LL_{x,y}^i
		\quad\Longrightarrow\quad j\le i.
	\end{equation}
	
	We prove the reverse inequality by induction
	on $\ell(x)$, with $y\in{}^I W$ fixed.
	For $x=e$, the assertion follows from the
	linearity of the ordinary Rouquier complex $F_y$.
	Suppose $x\ne e$, and assume that $B_aF_y$
	has a linear minimal model for every $a<x$.
	This induction is valid since $x\in W_I$
	implies $a\in W_I$ whenever $a<x$.
	Let
	\[
	K_x=\sigma_{\ge1}F_x
	\]
	be the brutal truncation of $F_x$ in positive
	degrees. By \cite[Theorem 6.9]{EW}, there is a
	termwise split exact sequence of complexes
	\begin{equation}\label{eq:truncation}
		0\longrightarrow K_x
		\longrightarrow F_x
		\longrightarrow B_x[0]
		\longrightarrow0,
	\end{equation}
	where
	\[
	K_x^k=
	\bigoplus_{a<x}B_a(k)^{\oplus n_{a,k}}
	\qquad(k\ge1).
	\]
	We claim that $K_xF_y$ is homotopy equivalent
	to a linear complex. Here $K_xF_y$ denotes
	the total complex of the tensor product
	$K_x\otimes_R F_y$, with terms
	\[
	(K_xF_y)^i=
	\bigoplus_{p+q=i}K_x^p\otimes_R F_y^q.
	\]

	For \(k\ge1\), put \(K_x^{\ge k}=\sigma_{\ge k}K_x\). There are exact sequences of complexes
	\[
	0\longrightarrow K_x^{\ge k+1}
	\longrightarrow K_x^{\ge k}
	\longrightarrow K_x^k[-k]
	\longrightarrow0.
	\]
	These sequences split in each cohomological degree. Tensoring with \(F_y\) and taking total complexes preserves their degreewise exactness. The resulting degreewise split exact sequences induce distinguished triangles in \(K^b(\mathcal B)\):
	\[
	K_x^{\ge k+1}F_y
	\longrightarrow K_x^{\ge k}F_y
	\longrightarrow K_x^k[-k]F_y
	\longrightarrow (K_x^{\ge k+1}F_y)[1].
	\]
	By the description of $K_x^k$, the third term
	is isomorphic to
	\[
	K_x^k[-k]F_y
	\cong
	\bigoplus_{a<x}
	\bigl((B_aF_y)(k)[-k]\bigr)^{\oplus n_{a,k}}.
	\]
	By induction, $B_aF_y$ has a linear minimal
	model $\LL_{a,y}$. Moreover,
	\[
	\bigl(\LL_{a,y}(k)[-k]\bigr)^i
	=\LL_{a,y}^{i-k}(k),
	\]
	whose indecomposable summands have internal
	shift $(i-k)+k=i$.
	Thus the third term of each triangle is
	homotopy equivalent to a linear complex.
	
	Since $K_x$ is bounded, we have
	$K_x^{\ge k}=0$ for sufficiently large $k$.
	Applying Lemma~\ref{lem:cones}(i) successively
	to the triangles above, we conclude that
	$K_xF_y$ is homotopy equivalent to a linear
	complex $A$. Now $\ell(xy)=\ell(x)+\ell(y)$, since
	$x\in W_I$ and $y\in{}^I W$.
	Tensoring \eqref{eq:truncation} with $F_y$
	therefore gives a distinguished triangle
	\[
	K_xF_y\longrightarrow F_{xy}
	\longrightarrow B_xF_y
	\longrightarrow (K_xF_y)[1],
	\]
	where we use $F_xF_y\simeq F_{xy}$.
	Consequently,
	\[
	B_xF_y\simeq\Cone(A\to F_{xy}).
	\]
	Both $A$ and $F_{xy}$ are linear, so
	Lemma~\ref{lem:cones}(ii) gives
	\[
	B_z(j)\summand\LL_{x,y}^i
	\quad\Longrightarrow\quad j\ge i.
	\]
	Together with \eqref{eq:upper}, this proves
	$j=i$, completing the induction.
	
	(ii) Let $y=s_1\cdots s_d$ be reduced,
	where $d=\ell(y)$.
	Each $F_{s_k}$ is supported in cohomological
	degrees $0$ and $1$. Hence
	\[
	B_xF_{s_1}\cdots F_{s_d}
	\]
	is supported in degrees $0,\ldots,d$.
	Its minimal model is obtained by removing
	contractible direct summands and therefore
	has the same degree bounds.
\end{proof}

\section{Signed positivity}
\label{sec:positivity}

\begin{proof}[Proof of Theorem~\ref{thm:positive}]
	(i) The KL basis of $\HH(W_I)$ agrees with the
	ambient KL basis. Since $y\in{}^I W$, for every
	$a\le x$ we have
	\[
	a\in W_I,\qquad
	\ell(ay)=\ell(a)+\ell(y),
	\]
	and $a<x$ implies $ay<xy$. Hence
	\[
	C_xH_y=H_{xy}+
	\sum_{a<x}h_{a,x}(v)H_{ay},
	\qquad h_{a,x}(v)\in v\Z[v].
	\]
	The transition from the standard basis to the
	KL basis is unitriangular, with off-diagonal
	coefficients in $v\Z[v]$. It follows that
	\begin{equation}\label{eq:triangular}
		\begin{aligned}
			r^I_{xy;x,y}&=1,\\
			r^I_{z;x,y}&=0
			&&\text{if }z\not\le xy,\\
			r^I_{z;x,y}&\in v\Z[v]
			&&\text{if }z<xy.
		\end{aligned}
	\end{equation}
	This proves (i).
	
	(ii) Since $\ch(F_s)=C_s-v=H_s$, we have
	$\ch(B_xF_y)=C_xH_y$.
	Taking Euler characters in
	Proposition~\ref{prop:linear} gives
	\begin{equation}\label{eq:euler}
		r^I_{z;x,y}(v)
		=\sum_{i=0}^{\ell(y)}(-1)^i m_{z,i}v^i.
	\end{equation}
	
	Consider the semilinear automorphism
	$\tau$ defined by
	\[
	\tau(v)=-v,\qquad
	\tau(H_w)=(-1)^{\ell(w)}H_w.
	\]
	It commutes with the bar involution, and the
	defining properties of the KL basis give
	$\tau(C_w)=(-1)^{\ell(w)}C_w$.
	Applying $\tau$ to \eqref{eq:inverse}, we obtain
	\begin{equation*}
		r^I_{z;x,y}(-v)
		=(-1)^{\ell(xy)-\ell(z)}
		r^I_{z;x,y}(v).
	\end{equation*}
	
	Comparing this identity with \eqref{eq:euler}, we see that
		\begin{equation}\label{eq:parity}
	m_{z,i}\ne0
	\quad\Longrightarrow\quad
	i\equiv\ell(xy)-\ell(z)\pmod2.
\end{equation}
	Consequently,
	\begin{equation}\label{eq:multiplicities}
		Q^I_{z;x,y}(v)
		=(-1)^{\ell(xy)-\ell(z)}
		r^I_{z;x,y}(v)
		=\sum_{i=0}^{\ell(y)}m_{z,i}v^i.
	\end{equation}
	Since $m_{z,i}\in\N$, this proves (ii).
	
	(iii) The degree bound follows from
	Proposition~\ref{prop:linear}(ii) and
	\eqref{eq:multiplicities}.
	The parity assertion follows from
	\eqref{eq:parity}.
\end{proof}

\section*{Disclosure of AI assistance}
The proofs in this paper were developed almost entirely by GPT-6 Astra.
The authors primarily reviewed and verified the arguments and their technical
details. The authors take full responsibility for the mathematical correctness
and content of this paper.

\end{document}